\documentclass[graybox]{svmult}

\usepackage{mathptmx}       
\usepackage{helvet}         
\usepackage{courier}        
\usepackage{graphicx}        
\usepackage{multicol}        
\usepackage[bottom]{footmisc}

\usepackage{cite}
\usepackage{fancyvrb}        

\usepackage{amsmath}
\usepackage{amsfonts}
\usepackage{arydshln}
\usepackage[caption=false]{subfig}

\DeclareMathOperator{\SchnaubeltCurl}{\mathrm{\mathbf{curl}}} 
\newcommand{\Schnaubeltdd}{\, \mathrm{d}}
\DeclareMathOperator{\SchnaubeltGrad}{\mathrm{\mathbf{grad}}} 

\newcommand{\SchnaubeltTM}[3]{$\text{TM}_{\scriptscriptstyle {#1}{#2}{#3}}$}
\newcommand{\SchnaubeltTE}[3]{$\text{TE}_{\scriptscriptstyle {#1}{#2}{#3}}$}

\begin{document}

\title*{A Comparison Between Different Formulations for Solving Axisymmetric Time-Harmonic Electromagnetic Wave Problems}
\titlerunning{Comparison of Different Formulations for Axisymmetric EM Wave Problems}
\author{Erik Schnaubelt 
\and
Nicolas Marsic 
\and
Herbert De Gersem 
}
\institute{Erik Schnaubelt, Nicolas Marsic, Herbert De Gersem \at Technische Universit\"at Darmstadt,
    Institut f\"ur Teilchenbeschleunigung und Elektromagnetische Felder (TEMF),
    64289 Darmstadt, Germany\\
\email{\{marsic,degersem\}@temf.tu-darmstadt.de}
  }

%
%
\maketitle

\abstract*{In many time-harmonic electromagnetic wave problems,
  the considered geometry exhibits an axial symmetry.
  In this case, by exploiting a Fourier expansion along the
  azimuthal direction, fully three-dimensional (3D) calculations
  can be carried out on a two-dimensional (2D) angular cross section
  of the problem, thus significantly reducing the computational effort.
 However, the transition from a full 3D problem to a 2D analysis introduces additional difficulties such as, among others, a singularity in the variational formulation. In this work, we compare and discuss different finite element formulations to deal with these obstacles. Particular attention is paid to spurious modes and to the convergence behavior when using high-order elements.}

\abstract{In many time-harmonic electromagnetic wave problems,
  the considered geometry exhibits an axial symmetry.
  In this case, by exploiting a Fourier expansion along the
  azimuthal direction, fully three-dimensional (3D) calculations
  can be carried out on a two-dimensional (2D) angular cross section
  of the problem, thus significantly reducing the computational effort.
 However, the transition from a full 3D problem to a 2D analysis introduces additional difficulties such as, among others, a singularity in the variational formulation. In this work, we compare and discuss different finite element formulations to deal with these obstacles. Particular attention is paid to spurious modes and to the convergence behavior when using high-order elements.}

\section{Introduction}
\label{schnaubelt:intro}

When treating a problem exhibiting axial symmetry,
a Fourier expansion along the azimuthal direction
can be exploited in order to restrict the computation
to a two-dimensional (2D) angular cross section of the geometry,
while still considering a fully three-dimensional (3D) solution~\cite{schnaubelt:Lacoste2000}. Therefore, these methods are also referred to as \textit{quasi-3D} or \textit{2.5D} methods.
Let us consider a cylindrical coordinate system $(r, \varphi, z)$,
and let us expand the electric field $\vec{e}(r, \varphi, z)$
into a Fourier series along~$\varphi$:
\begin{equation*}
\def\arraystretch{1.2}
  \vec{e}(r, \varphi, z) =
  \left[
    \begin{array}{c}
      e^0_r(r,z)\\
      e^0_\varphi(r,z)\\
      e^0_z(r,z)
    \end{array}
  \right]
  +
  \sum_{m=1}^\infty
  \left(
    \left[
      \begin{array}{c}
        e^m_r(r,z)      \cos(m\varphi)\\
        e^m_\varphi(r,z)\sin(m\varphi) \\
        e^m_z  (r,z)    \cos(m\varphi)
      \end{array}
    \right]%
    \hspace{-3pt}+\hspace{-3pt}%
    \left[
      \begin{array}{c}
        e^{-m}_r (r,z)     \sin(m\varphi)\\
        e^{-m}_\varphi(r,z)\cos(m\varphi) \\
        e^{-m}_z    (r,z)  \sin(m \varphi)
      \end{array}
    \right]
  \right),
\end{equation*}
where the Fourier coefficients
$\vec{e}^n(r,z) = \big[e^n_r, e^n_\varphi,e^n_z\big]^T$ with $n\in \mathbb{Z}$
are functions of the radial and axial coordinates \emph{only}.
Furthermore, by exploiting the orthogonality of the trigonometric functions,
we can write the Maxwell eigenvalue problem
for an axisymmetric cavity $V$ with perfect electric conducting boundaries
as~\cite{schnaubelt:Lacoste2000}:
\begin{equation}
  \label{schnaubelt:eig}%
  \left\{
    \begin{array}{l}
      \text{For a given mode $n \in \mathbb{Z}$, find the eigenpairs}~(\vec{e}^n, \omega^2)~%
      \text{with}~\vec{e}^n\in\mathcal{S}^n(\Omega):
      \\
      \displaystyle\int_\Omega
      \vec{\mu_r}^{-1}
      \SchnaubeltCurl_n\vec{e}^n\cdot\SchnaubeltCurl_n\vec{e}^{n^\prime} \Schnaubeltdd{}\Omega
       -
      \frac{\omega^2}{c_0^2}
      \int_\Omega \vec{\varepsilon_r} \, \vec{e}^n\cdot\vec{e}^{n^\prime}\Schnaubeltdd{}\Omega
      = 0 \quad
      \forall\vec{e}^{n^\prime}\hspace{-3pt}\in\mathcal{S}^n(\Omega),
    \end{array}
  \right.
\end{equation}
with $\vec{\varepsilon_r} = \vec{\varepsilon_r}(r,z)$ and
$\vec{\mu_r} = \vec{\mu_r}(r,z)$ the (possibly anisotropic) relative electric permittivity and magnetic permeability tensors of the medium, $\Omega$ a 2D angular cross section of $V$, $\Schnaubeltdd{}\Omega = r\Schnaubeltdd{}r\Schnaubeltdd{}z$,
$\omega$~the angular frequency, $c_0$ the speed of light in vacuum,
$\mathcal{S}^n(\Omega)$ the function space of the $n$\textsuperscript{th}
Fourier coefficient and
\begin{equation*}
  \SchnaubeltCurl_n\vec{e}^n =
  \left[
    \begin{array}{c}
      -r^{-1}\big(ne_z^n+\partial_z(re_\varphi^n)\big)
      \\[3pt]
      \partial_z{}e_r^n - \partial_r{}e_z^n
      \\[3pt]
      +r^{-1}\big(ne_r^n+\partial_r(re_\varphi^n)\big)
    \end{array}
  \right].
\label{schnaubelt:curl}
\end{equation*}

\section{Well-Posed Variational Formulation}
In order to construct an appropriate subspace of $\mathcal{S}^n(\Omega)$
and in order to account for the singular behavior of $\SchnaubeltCurl_n$ at ${r=0}$,
two strategies have been proposed in the literature.

\subsection{Non-Classical Conditions Along the Symmetry Axis}
\label{schnaubelt:lee}

A first approach consists in taking the unknown fields $e_\varphi^{\star,n}=re_\varphi^n\in{}H^1(\Omega)$ and  $\vec{e}_{rz}^n = \big[e_r^n, e_z^n\big]^T\in{}H(\SchnaubeltCurl,\Omega)$~\cite{schnaubelt:Lee1993}
together with \emph{non-classical discrete conditions at the symmetry axis}~\cite{schnaubelt:Chinellato2005b}.
By following this strategy, all integrals are well-posed but exhibit singular integrands, hence requiring either \textit{i)} a classical Gau\ss{} quadrature with a large number of quadrature points or \textit{ii)} specialized quadrature rules~\cite{schnaubelt:Chinellato2005b} which differ from element to element, thus preventing the use of fast assembly techniques~\cite{schnaubelt:Marsic2015}.
In what follows, this approach will be further referred to as transformation ``TA''.

\subsection{Direct Construction of a Subspace of $\mathcal{S}^n(\Omega)$}

Another approach consists in \emph{directly constructing an appropriate subspace of $\mathcal{S}^n(\Omega)$} such that the variational formulation is guaranteed to be always well-posed,
as shown in~\cite{schnaubelt:Lacoste2000, schnaubelt:Dunn2006, schnaubelt:Cambon2012, schnaubelt:Oh2015} for instance, thus avoiding the need of non-classical conditions on the symmetry axis. To this end, the unknowns $e_\varphi^n$ and $\vec{e}_{rz}^n$ are transformed into $\tilde{u}^n\in{}H^1(\Omega)$ and $\vec{\tilde{U}}^n\in{}H(\SchnaubeltCurl,\Omega)$ by following the methodology shown in Table~\ref{schnaubelt:tabtransform}, with  $\SchnaubeltGrad_{rz} e_\varphi^n = \big[\partial_r{}e_\varphi^n, \partial_z{}e_\varphi^n\big]^T$ and $\vec{\hat{r}}$ the unit vector along the $r$-axis.

\begin{table}[htb]
  \def\arraystretch{2}
  \setlength\tabcolsep{.25cm}
  \centering
  \caption{Different transformations for constructing a subspace of $\mathcal{S}^n(\Omega)$.}
  \label{schnaubelt:tabtransform}
  \begin{tabular}{llll}
    \svhline
    \textbf{Mode} & \textbf{Transf. TB}~\cite{schnaubelt:Oh2015} & \textbf{Transf. TC}$(\alpha,\beta)$~\cite{schnaubelt:Cambon2012} & \textbf{Transf. TD}~\cite{schnaubelt:Dunn2006} \\
    \hline
    $n = 0$    & $\tilde{u}^0=e_\varphi^0$
               & $r^\beta\tilde{u}^0=re_\varphi^0$
               & $\tilde{u}^0=e_\varphi^0$\\
               & $\vec{\tilde{U}}^0=\vec{e}_{rz}^0$
               & $\vec{\tilde{U}}^0=\vec{e}_{rz}^0$
               & $\vec{\tilde{U}}^0=\vec{e}_{rz}^0$\\
    \hdashline
    $n=\pm{}1$ & $\tilde{u}^{\pm{}1}=e_\varphi^{\pm{}1}$
               & $r^\beta\tilde{u}^{\pm{}1}=re_\varphi^{\pm{}1}$
               & $\tilde{u}^{\pm{}1}=e_\varphi^{\pm{}1}$\\
               & $\vec{\tilde{U}}^{\pm{}1}=\displaystyle\frac{n}{r}\vec{e}_{rz}^{\pm{}1}+\frac{e_\varphi^{\pm{}1}}{r}\vec{\hat{r}}$
               & $r^\alpha\vec{\tilde{U}}^{\pm{}1}=\pm{}\vec{e}_{rz}^{\pm{}1} + \SchnaubeltGrad_{rz}(re_\varphi^{\pm{}1})$
               & $\vec{\tilde{U}}^{\pm{}1}=\displaystyle\frac{n}{r}\vec{e}_{rz}^{\pm{}1}+\frac{e_\varphi^{\pm{}1}}{r}\vec{\hat{r}}$\\[3pt]
    \hdashline
    $|n|>1$    & $\tilde{u}^n=e_\varphi^n$
               & $r^\beta\tilde{u}^n=re_\varphi^n$
               & $\tilde{u}^n=e_\varphi^n$\\
               & $\vec{\tilde{U}}^n=\displaystyle\frac{n}{r}\vec{e}_{rz}^n+\frac{e_\varphi^n}{r}\vec{\hat{r}}$
               & $r^\alpha\vec{\tilde{U}}^n=n\vec{e}_{rz}^n + \SchnaubeltGrad_{rz}(re_\varphi^n)$
               & $\vec{\tilde{U}}^n=\displaystyle\frac{n}{r}\vec{e}_{rz}^n$\\[3pt]
    \svhline
  \end{tabular}
\end{table}

The parameters $\alpha$ and $\beta$ in TC$(\alpha,\beta)$ must satisfy, according to~\cite{schnaubelt:Cambon2012}, the constraints shown in Table~\ref{schnaubelt:tabalphabeta}.
Furthermore, some transformations need an additional \emph{homogeneous Dirichlet condition} at $r = 0$,
as shown in Table~\ref{schnaubelt:tabdirichlet}.
Finally, for appropriate choices of $\alpha$ and $\beta$, TC$(\alpha,\beta)$ leads to \emph{polynomial} integrands (see Section~\ref{schnaubelt:sceconvergence} for more details). This property is also met by TB for $n\neq{}0$ and TD for $n=\pm{}1$.

\begin{table}[htb]
  \def\arraystretch{2}
  \setlength\tabcolsep{.1cm}
  \centering
  \caption{Constraints on $\alpha$ and $\beta$ for TC$(\alpha, \beta)$ according to~\cite{schnaubelt:Cambon2012}.}
  \label{schnaubelt:tabalphabeta}
  \begin{tabular}{p{3cm}p{3cm}p{3cm}}
    \svhline
    $ n = 0$ & $n = \pm 1$ & $\lvert n \rvert > 1$ \\
    \hline
    $\beta \geq 0.5$ & $\alpha \geq 0.5$ and $\beta = 1$ & $\alpha \geq 0.5$ and $\beta > 0$ \\
    \svhline
  \end{tabular}
\end{table}
\begin{table}[htb]
  \def\arraystretch{2}
  \setlength\tabcolsep{.1cm}
  \centering
  \caption{Conditions on the symmetry axis.}
  \label{schnaubelt:tabdirichlet}
  \begin{tabular}{llll}
    \svhline
    \textbf{Mode} & \textbf{Transf. TB}~\cite{schnaubelt:Oh2015} & \textbf{Transf. TC}$(\alpha,\beta)$~\cite{schnaubelt:Cambon2012} & \textbf{Transf. TD}~\cite{schnaubelt:Dunn2006} \\
    \hline
    $n = 0$    & $\tilde{u}^0=0$
               & $\tilde{u}^0=0$ if $\beta\in[0.5, 1.5[$, none otherwise
               & $\tilde{u}^0=0$\\
    $n=\pm{}1$ & None
               & None
               & None\\
    $|n|>1$    & $\tilde{u}^n=0$
               & $\tilde{u}^n=0$ if $\beta\in]0, 1]$, none otherwise
               & $\tilde{u}^n=0$\\
    \svhline
  \end{tabular}
\end{table}

\section{Comparison and Discussion of the Quasi-3D Methods}

As already stated, this work compares the aforementioned transformations to treat the eigenvalue problem \eqref{schnaubelt:eig}. To this end, they were implemented in a homemade high-order finite element (FE) code\footnote{See \textit{https://gitlab.onelab.info/gmsh/small\_fem/blob/master/simulation/Quasi3D.cpp}.}. All following numerical experiments are performed on a pillbox cavity with $\vec{\mu_r} = \vec{\varepsilon_r} = 1$ for which closed-form solutions are well-known \cite{schnaubelt:Wangler2008}.

\subsection{Spurious Modes and High-Order FE Discretizations}

Let us start our comparison by determining if the methods discussed previously can avoid spurious modes.
As we search the azimuthal unknown ($e_\varphi^{*,n}$ for TA and $\tilde{u}^n$ for TB, TC and TD) in a \emph{finite} subspace of $H^1(\Omega)$ of polynomial order $q$ and the in-plane unknown ($\vec{e}_{rz}$ for TA and $\vec{\tilde{u}}^n$ for TB, TC and TD) in a \emph{finite} subspace of $H(\SchnaubeltCurl,\Omega)$ of polynomial order\footnote{In this paper, we consider subspaces of $H(\SchnaubeltCurl,\Omega)$ with both irrotational and rotational functions.} $p$,
the dimension of each subspace must be selected with care.
In particular, in order to satisfy the \emph{exactness of the discrete de Rham sequence} \cite{schnaubelt:simona2019isogeometric}, one must impose that $q = p+1$~\cite{schnaubelt:Demkowicz2000}.

\begin{figure}[ht]
   \centering
     \includegraphics[]{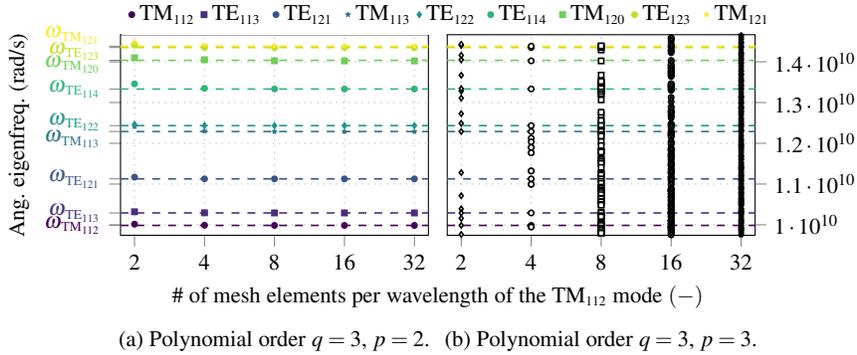}
   \caption{Part of the spectrum of a pillbox cavity obtained with TB and $n=1$.}
   \label{schnaubelt:spurMod}%
\end{figure}

In order to validate this choice, we ran multiple numerical tests with the different transformations, different modes $n$ and different values for $p$ and $q$.
As a result, we observed that, \emph{apart from TD}, all eigenspectra were free of spurious modes when $q=p+1$.
Interestingly, we also observed no spurious modes when $q>p+1$.
On the other hand, spurious modes were systematically observed when $q<p+1$, and when transformation TD was used with $|n|>1$ (for all possible values of $p$ and $q$).
For this reason, TD will not be investigated further.
As an illustration, Fig.~\ref{schnaubelt:spurMod} shows a part of the numerical spectrum of a pillbox cavity for $n = 1$ and different mesh densities. It was computed with TB, once for $q = 3, p = 2$ and once for $q = p = 2$.

When $n = 0$, the in-plane and azimuthal unknowns are \emph{decoupled} from each other~\cite{schnaubelt:Cambon2012}.
Therefore, $q$ and $p$ can be chosen \emph{independently}.
In the following, the symbol ``$\ast$'' will be used to indicate a result independent of $p$ or $q$ (see Fig.~\ref{schnaubelt:alphabeta:n0}).

\subsection{Convergence Results for Higher Order Finite Elements}
\label{schnaubelt:sceconvergence}

\begin{figure}[hbt]
   \centering
      \includegraphics[]{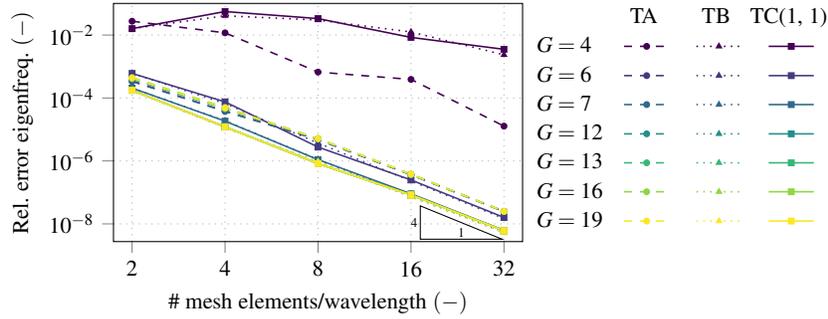}
   \caption{Convergence results when computing the eigenfrequency of the \SchnaubeltTE{1}{1}{1} mode of a pillbox cavity with TA, TB and TC(1, 1), together with $q = 3$, $p = 2$ and $G$ quadrature points.}
   \label{schnaubelt:highOrd}
\end{figure}

In this subsection, the convergence behavior of the  different formulations in combination with high-order basis functions will be compared. As an example, the evolution of the relative error between the numerically computed eigenfrequency and its analytical counterpart is shown in Fig.~\ref{schnaubelt:highOrd} for different mesh densities and different number of quadrature points~$G$. A second-order FE method with $q = 3, p = 2$ is used, hence resulting in an expected convergence slope of 4 \cite{schnaubelt:Sauter2010}. This slope is indeed achieved for TA, TB and TC(1, 1) once the number of quadrature points passes a certain threshold.
Again, this is not an isolated case but can be observed for all $n$ and for different element orders.

Let us also stress that TB and TC(1, 1) \textit{i)} yield a lower relative error than TA and \textit{ii)} depend less on $G$ than TA.
This last observation can be easily explained: as TB and TC(1, 1) lead to \emph{polynomial integrands}, the final solution is independent of $G$, at least for a $G$ sufficiently large to integrate a polynomial of the given order \emph{exactly}.

The parameters of TC($\alpha, \beta$) need to meet the following criteria to yield polynomial integrands in the variational formulation~\eqref{schnaubelt:eig}: \textit{i)} $\alpha$ and $\beta$ must be \textit{multiples} of $0.5$, \textit{ii)} their sum must be an \textit{integer} and \textit{iii)} $\beta \geq 1.5$ for $n = 0$ and $\alpha \geq 0.5, \beta \geq 0.5$ for $n \neq 0$. 

\subsection{Influence of $\alpha$ and $\beta$ on the Convergence Behavior}
\label{schnaubelt:secalphabeta}

\begin{figure}[tbh]
\subfloat[][\SchnaubeltTE{0}{2}{2} mode with $q = 4$, $p = \ast$.]{ 
 \centering
        \includegraphics[]{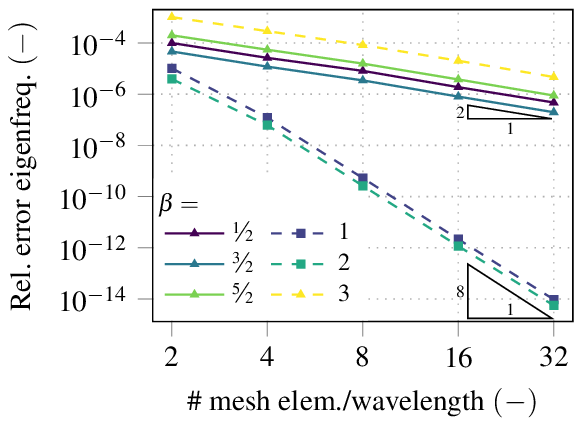}
    \label{schnaubelt:alphabeta:n0}%
}
\hspace{-5pt}
\subfloat[][\SchnaubeltTM{1}{1}{1} mode with $q = 4$, $p = 3$.]{ 
 \centering
         \includegraphics[]{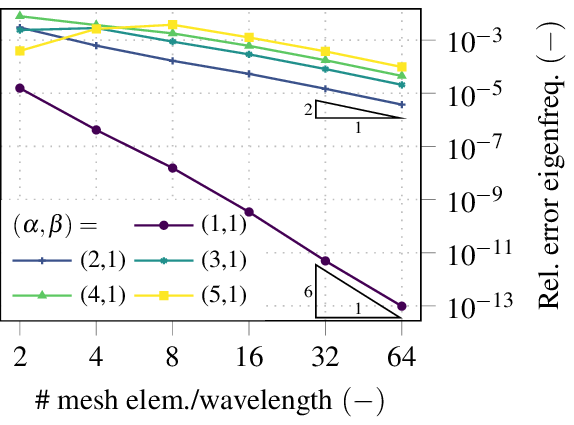}
}
\caption{Convergence rate of TC$(\alpha, \beta)$ for different values of $\alpha$ and $\beta$, as allowed by Table~\ref{schnaubelt:tabalphabeta}.}
\label{schnaubelt:alphabeta}%
\end{figure}

Let us now focus on the transformation TC$(\alpha, \beta)$, and let us carry out a convergence test similar to the previous section.
However, now, the influence of the parameters $\alpha$ and $\beta$ (chosen according to Table~\ref{schnaubelt:tabalphabeta}) on the convergence rate will be investigated.
The results of this numerical experiment are displayed in Fig.~\ref{schnaubelt:alphabeta}. 
As it can be observed directly, while all choices converge towards the sought eigenvalue, only particular pairs $(\alpha, \beta)$ exhibit the expected convergence rate.
This behavior has been observed for other choices of $(n,p,q)$ with $q=p+1$ as well. 

This behavior can be easily explained if we assume that $\vec{e}^n\in{}C^\infty$ in the vicinity of the symmetry axis.
This assumption is of course restrictive, but applies to the pillbox cavity~\cite{schnaubelt:Wangler2008}, and gives already a good insight into the underlying numerical mechanisms.
In what follows, only the case $n=\pm{}1$ will be discussed, but the same methodology applies to the other cases.

Let us start by expanding $\vec{e}^{\pm{}1}$ into a Taylor series in the vicinity of $r=0$ and $z=z_0$. As $e_z^{\pm{}1}=0$ at $r=0$ (see~\cite{schnaubelt:Lacoste2000}), we have:
\begin{equation*}
  \left\{
    \begin{array}{rcl}
      e_\varphi^{\pm{}1}(r,z) & = & a_0 + a_1^rr+a_1^z(z-z_0)+a_2^{rr}r^2+a_2^{zz}(z-z_0)^2+2a_2^{rz}r(z-z_0)+\dots,\\
      e_z^{\pm{}1}(r,z)       & = & b_1^rr+b_2^{rr}r^2+\dots,\\
      e_r^{\pm{}1}(r,z)       & = & c_0 + c_1^rr+c_1^z(z-z_0)+c_2^{rr}r^2+c_2^{zz}(z-z_0)^2+2c_2^{rz}r(z-z_0)+\dots \, .
    \end{array}
  \right.
\end{equation*}
Then, by exploiting the definition of $\vec{\tilde{U}}^{\pm{}1}$ (see Table~\ref{schnaubelt:tabtransform}), we can write:
\begin{equation*}
  \left\{
    \begin{array}{lcl@{}c@{}l}
      r^\alpha \tilde{U}_r^{\pm{}1} & = & \pm{}e_r^{\pm{}1} + & e_\varphi^{\pm{}1} & + r(a_1^r+2a_2^{rr}r+2a_2^{rz}(z-z_0)+\dots),\\
      r^\alpha \tilde{U}_z^{\pm{}1} & = & \pm{}e_z^{\pm{}1} + & 0           & + r(a_1^z+2a_2^{zz}(z-z_0)+2a_2^{rz}r+\dots).\\
    \end{array}
  \right.
\end{equation*}
Afterwards, by taking the limit $r\to{}0$, and by exploiting the conditions $e_r^{\pm{}1}\pm{}e_\varphi^{\pm{}1}=0$ at $r=0$ and $e_z^{\pm{}1}=0$ at $r=0$ (see~\cite{schnaubelt:Lacoste2000}), it follows that:
\begin{equation*}
  \left\{
    \begin{array}{lclcl}
      \displaystyle\lim_{r\to{}0}\tilde{U}_r^{\pm{}1} & = & \displaystyle\lim_{r\to{}0}r^{-\alpha}\Big[0+r(a_1^r+2a_2^{rr}r+2a_2^{rz}(z-z_0)+\dots)\Big] & = & \displaystyle\lim_{r\to{}0}r^{1-\alpha}f(z),\\
      \\[-6pt]
      \displaystyle\lim_{r\to{}0}\tilde{U}_z^{\pm{}1} & = & \displaystyle\lim_{r\to{}0}r^{-\alpha}\Big[0+r(a_1^z+2a_2^{zz}(z-z_0)+2a_2^{rz}r+\dots)\Big]       & = & \displaystyle\lim_{r\to{}0}r^{1-\alpha}g(z),
    \end{array}
  \right.
\end{equation*}
where $f$ and $g$ are functions of $z$.
In other words, the unknown $\vec{\tilde{U}}^{\pm{}1}$ behaves as $r^{1-\alpha}$ as $r\to{}0$.
Therefore, for the sought FE solution to be differentiable at $r=0$, and because of the constraint enforced by Table~\ref{schnaubelt:tabalphabeta} ($\alpha\geq{}0.5$), \emph{we need to impose that $\alpha~=~1$}.
Any other choice (in accordance with Table~\ref{schnaubelt:tabalphabeta}) will lead to a non-differentiable $\vec{\tilde{U}}^n$, jeopardizing thus the convergence of the FE scheme.
By applying the same strategy to the cases $n=0$ and $|n|>1$, and by taking into account the restrictions imposed in Table~\ref{schnaubelt:tabalphabeta}, the set of allowed couples $(\alpha,\beta)$ must be further narrowed, as shown in Table~\ref{schnaubelt:tabalphabetahp}.
To the best of our knowledge, these last results have been derived for the first time.

\begin{table}[htb]
  \def\arraystretch{2}
  \setlength\tabcolsep{.1cm}
  \centering
  \caption{Values of $\alpha$ and $\beta$ for TC$(\alpha, \beta)$ leading to a high FE convergence rate.}
  \label{schnaubelt:tabalphabetahp}
  \begin{tabular}{p{3cm}p{3cm}p{3cm}}
    \svhline
    $ n = 0$ & $n = \pm 1$ & $\lvert n \rvert > 1$ \\
    \hline
    $\beta\in\{1,2\}$ & $\alpha=1$ and $\beta = 1$ & $\alpha\in\{1,2\}$ and $\beta\in\{1,2\}$ \\
    \svhline
  \end{tabular}
\end{table}

For the integrands of the variational formulation~\eqref{schnaubelt:eig} to be polynomial,
we further need to impose the following restriction: $\beta=2$ when $n=0$ (see Section \ref{schnaubelt:sceconvergence}).
Moreover, as the number of quadrature points depends on the order of the integrands, it is preferable to select the smallest acceptable $(\alpha, \beta)$ couple.
Therefore, the values given in Table~\ref{schnaubelt:tabalphabetarecommand} are recommended.
Let us finally note that, apart from the case $n=0$, these recommendations are in accordance with~\cite{schnaubelt:Cambon2012}, where the best $(\alpha, \beta)$
couples were determined on the basis of numerical experiments.

Concerning the case $n=0$, the choice $\beta=1$ leads to slightly more accurate results for a given mesh density in a numerical experiment carried out in~\cite{schnaubelt:Cambon2012}. However, our numerical experiments do not confirm these results as the choice $\beta = 2$ leads to a slightly lower relative error for the same mesh density (see Fig. \ref{schnaubelt:alphabeta:n0}). This behavior is supported by the fact that the integrands are polynomial for the latter choice $\beta = 2$.  

\begin{table}[htb]
  \def\arraystretch{2}
  \setlength\tabcolsep{.1cm}
  \centering
  \caption{Recommended choice of $\alpha$ and $\beta$ for TC$(\alpha, \beta)$.}
  \label{schnaubelt:tabalphabetarecommand}
  \begin{tabular}{p{3cm}p{3cm}p{3cm}}
    \svhline
    $ n = 0$ & $n = \pm 1$ & $\lvert n \rvert > 1$ \\
    \hline
    $\beta=2$ & $\alpha=1$ and $\beta = 1$ & $\alpha=1$ and $\beta=1$ \\
    \svhline
  \end{tabular}
\end{table}

\section{Conclusion} 

This paper compared four different transformations to treat three-dimensional time-har\-mo\-nic electromagnetic wave problems in axisymmetric geometries proposed in the literature.
We first determined numerically that the transformations TA, TB and TC$(\alpha, \beta)$ lead to eigenvalue problems which are free of spurious modes, while the transformation TD exhibits spurious modes when $|n|>1$.
We then compared numerically the accuracy of TA, TB and TC$(\alpha, \beta)$, and found that TB and TC$(\alpha, \beta)$ produce the most accurate results for a given mesh density.
Finally, we analyzed theoretically the convergence rate of TC$(\alpha, \beta)$ for different values of $\alpha$ and $\beta$ in a high-order FE context, and determined new restrictions on $\alpha$ and $\beta$.

\begin{acknowledgement}
The authors would like to express their gratitude to Abele Simona for his valuable advice and the fruitful discussions on axisymmetric problems. 
\end{acknowledgement}

\bibliographystyle{spmpsci}
\bibliography{schnaubelt_scee.bib}

\end{document}